\documentclass[1p]{elsarticle}
\usepackage{amsfonts}
\usepackage{amssymb}

\def\NN{{\mathbb N}}
\def\RR{{{\mathbb R}}}

\def\SSS{{\mathbb S}}

\def\re{\mathop{\rm Re}\, }

\let\a=\alpha
\let\b=\beta

\let\O=\Omega

\def\cA{{\mathcal A}}

\def\cF{{\mathcal F}}

\def\cS{{\mathcal S}}

\def\pa{\partial}

\newtheorem{lemm}{Lemma} [section]
\newtheorem{theo}{Theorem} [section]
\newtheorem{defi}{Definition}[section]

\newtheorem{prop}{Proposition}[section]
\newtheorem{rema}{Remark}[section]

\begin{document}

\title{Ultra-analytic effect of Cauchy problem\\
for a class of kinetic equations}
\author{Yoshinori MORIMOTO}
\ead{morimoto\@@math.h.kyoto-u.ac.jp}
\address{Graduate School of Human and Environmental Studies\\
Kyoto University, Kyoto, 606-8501, Japan}
\author{Chao-Jiang XU}
\ead{Chao-Jiang.Xu\@@univ-rouen.fr}
\address{Universit\'e de Rouen, UMR-6085, Math\'ematiques\\
76801 Saint Etienne du Rouvray, France\\
and\\
School of Mathematics and Statistics, Wuhan University, Wuhan
430072, China }
\maketitle

\begin{center}
{\bf Abstract}
\end{center}

The smoothing effect of the  Cauchy problem for a class of kinetic
equations is studied. We firstly consider the spatially homogeneous
nonlinear Landau equation with Maxwellian molecules and
inhomogeneous linear Fokker-Planck equation to show the
ultra-analytic effects of the Cauchy problem. Those smoothing effect
results are optimal and similar to heat equation. In the second
part, we study a model of spatially inhomogeneous linear Landau
equation with Maxwellian molecules, and show the analytic effect of
the Cauchy problem.

\bigbreak \noindent {Key words:} {Landau equation,  Fokker-Planck
equation,  ultra-analytic effect of Cauchy problem.}

\bigbreak \noindent{ AMS Classification: 35A05, 35B65, 35D10, 35H20,
76P05, 84C40}

\section{Introduction}\label{s1}

It is well known that the Cauchy problem of heat equation possesses the
ultra-analytic effect phenomenon, namely, if $u(t, x)$ is the
solution of the following Cauchy problem :
$$
\left\{\begin{array}{ll} \partial_t u-\Delta_x u=0,&x\in\RR^d;
\,\,\,\,t>0\\
u|_{t=0}=u_0\in L^2(\RR^d),&
\end{array}\right.
$$
then under the uniqueness hypothesis, the solution $u(t,\,
\cdot\,)=e^{t\Delta_x}u_0$ is an ultra-analytic function for any
$t>0$. We give now the definition of function spaces $\cA^{s}(\O)$
 where $\O$ is an open subset of $\RR^d$.

\begin{defi}\label{defi0.1}
For $0<s<+\infty$, we say that $f\in \cA^s(\O)$, if $f\in
C^\infty(\O)$, and there exists $C>0, N_0>0$ such that
$$
\|\partial^\a f\|_{L^2(\O)}\leq
C^{|\a|+1}(\a!)^s,\,\,\,\,\,\,\,\,\,\,\forall\,\,\a\in\NN^d,\,\,\,|\a|\geq
N_0.
$$
If the boundary of $\Omega$ is smooth, by using Sobolev embedding
theorem, we have the same type estimate with $L^2$ norm replaced by any $L^p$ norm for
$2 < 
p \leq +\infty$.  On the whole space $\O=\RR^d$, it is also
equivalent to
$$
e^{c_0(-\Delta)^{\frac 1{2s}}} (\partial^{\b_0} f)\in L^2(\RR^d)
$$
for some $c_0>0$ and $ \b_0\in\NN^d$, where $
e^{c_0(-\Delta)^{\frac 1{2s}}}$ is the Fourier multiplier defined by
$$
e^{c_0(-\Delta)^{\frac 1{2s}}} u(x)=\cF^{-1}\Big(
e^{c_0|\xi|^{\frac 1{s}}} \hat{u}(\xi)\Big).
$$
\end{defi}

If $s=1$, it is usual analytic function. If $s>1$, it is Gevrey
class  function. For $0<s<1$, it is called ultra-analytic function.
Notice that all polynomial functions are ultra-analytic for any $s >0$.

It is obvious that  if $u_0\in L^2(\RR^d)$ then,  for any $t>0$ and
any $k\in \NN$, we have  $u(t, \, \cdot\,)=e^{-t(-\Delta_x)^k}u_0
\in\cA^{\frac 1{2k}}(\RR^d)$, namely, there exists $C>0$ such that
for any $m\in\NN,$
\begin{eqnarray*}
\|(t^m \partial_x^{2k m}) u(t, \,\cdot\,)\|_{L^2(\RR^d)} &\leq&
C^{k\, m}\|(t (-\Delta_x)^k)^m u(t, \,\cdot\,)\|_{L^2(\RR^d)}\\
&\leq& \|u_0\|_{L^2(\RR^d)} C^{k\, m} m\, !\leq \tilde{C}^{2k\, m+1}
\big((2km)\, !\big)^{\frac{1}{2k}}\,,
\end{eqnarray*}
where $\partial^{2k m}_x=\sum_{|\alpha|=2 k m,
\alpha\in\NN^d}\partial_x^\alpha$. We say that the diffusion
operators $(-\Delta_x)^k$ possess the ultra-analytic effect property
if $k> 1/2$, the analytic effect property if $k=1/2$ and  the Gevrey
effect property if $0<k<1/2$.

We study the Cauchy problem for spatially homogeneous Landau
equation
\begin{equation}\label{H-landau}\left\{\begin{array}{l}
f_t=Q(f,\, f)\equiv\nabla_v\big(\bar{a}(f)\cdot\nabla_v f-\bar{b}(f)
f\big),
\,\,\,\,\,\,\,\,v\in\RR^d,\,\, t>0,\\
f|_{t=0}=f_0 {  ,} \end{array}\right.
\end{equation}
where $\bar{a}(f) =(\bar{a}_{i j}(f))$ and
$\bar{b}(f)=(\bar{b}_1(f),\, \cdots\,,\, \bar{b}_d(f)\,)$ are
defined as follows (convolution is w.\,r.\,t.\, the variable
$v\in\RR^d$)
\[
\bar{a}_{i j}(f)=a_{i j}\star f,\,\,\,\,\,\, \bar{b}_j(f)=\sum^{{
d}}_{i=1} \big(\partial_{v_i}a_{i j}\big)\star f\, , \,\,\,\,\,\, i,
j=1,\, \cdots,\, d,
\]
with
\[
a_{i j}(v)=\left(\delta_{i j}-\frac{v_i
v_j}{|v|^2}\right)|v|^{\gamma+2},\,\,\,\,\,\,\,\,\,\gamma\in[-3, 1].
\]
We consider hereafter only the Maxwellian molecule case which corresponds
to $\gamma=0$.
We introduce also the notation, for $l\in\RR,\, L^p_l(\RR^d)=
\{f;\, (1+|v|^2)^{l/2} f\in L^p(\RR^d)\}$ is the weighted function space.

We prove the following ultra-analytic effect results for the nonlinear 
Cauchy problem (\ref{H-landau}).
\begin{theo}\label{theo0.1}
Let $f_0\in L^2(\RR^d)\bigcap L^1_2(\RR^d)$ and $0<T\leq +\infty$.
If $f(t, x)>0$  and $f\in L^\infty(]0,\, T[;\, L^2(\RR^d)\bigcap
L^1_2(\RR^d))$ is a weak solution of the Cauchy problem
(\ref{H-landau}), then for any $0<t<T$, we have
$$
f(t, \,\cdot\,)\in \cA^{1/2}(\RR^d),
$$
and moreover, for any $0<T_0<T$, there exists $c_0>0$ such that for
any $0<t\leq T_0$
\begin{equation}\label{estimtion1}
\left\|e^{-c_0t \triangle_v} {f}
(t,\,\cdot\,)\right\|_{L^2(\RR^d)}\leq e^{\frac d 2
t}\|f_0\|_{L^2(\RR^d)}.
\end{equation}
\end{theo}

In \cite{MUXY-DCDS}, they proved the Gevrey regularity effect of {
the} Cauchy problem for linear spatially homogeneous non-cut-off
Boltzmann equation. By a careful revision for the proof of Theorem
1.2 of \cite{MUXY-DCDS}, one can also prove that the solution of {
the} Cauchy problem (1.10) in \cite{MUXY-DCDS} belongs to
$\cA^{\frac 1{2\alpha}}(\RR^d)$ for any $t>0$, where $0<\alpha<1$ is
the order of singularity of collision kernel of Boltzmann operator.
Hence, if $\alpha\geq 1/2$, there is also the ultra-analytic effect
phenomenon. Now the above Theorem \ref{theo0.1} shows that, for
Landau equation, the ultra-analytic effect phenomenon holds in nonlinear case, 
which is an optimal regularity {  result.}

\smallbreak The ultra-analytic effect property is also true for
the Cauchy problem of the following generalized Kolmogorov operators
$$
\left\{\begin{array}{ll} \partial_t u+v\,\cdot\,\nabla_x u
+\left(-\Delta_v\right)^\alpha
u=0,&(x, v)\in\RR^{2d};\,\,\,\,
t>0\\
u|_{t=0}=u_0\in L^2(\RR^{2d}),&
\end{array}\right.
$$
where $0<\alpha<\infty$, and the classical Kolmogorov operators is
corresponding to $\alpha=1$.
By Fourier transformation, the explicit solution of the above Cauchy
problem is given by
$$
\hat{u}(t, \eta, \xi)=e^{-\int^t_0|\xi+s\eta|^{2\alpha} ds }\hat{u}_0(\eta,
\xi+t\eta).
$$
Since there exists $c_\alpha>0$ (see Lemma \ref{lemm4.1} below) such that
\begin{equation}\label{0.1}
c_\alpha \,(t|\xi|^{2\alpha}+t^{2\alpha+1}|\eta|^{2\alpha})\leq
\int^t_0|\xi+s\eta|^{2\alpha}ds,
\end{equation}
we have
$$
e^{c_\alpha(t(-\Delta_v)^\alpha+t^{2\alpha+1}(-\Delta_x)^{{
\alpha}})} u(t, \,\cdot,\,\cdot\,)\in L^2(\RR^{2d}),
$$
i. e. $u(t,\,\cdot,\,\cdot\,)\in \cA^{1/(2\alpha)}(\RR^{2d})$ for any $t>0$.

Notice that this ultra-analytic (if $\alpha>1/2$) effect phenomenon
is similar to heat equations of $(x, v)$ variables. That is, this
means $\,\,v\,\cdot\,\nabla_x +(-\Delta_v)^\alpha\,$ is equivalent
to $(-\Delta_x)^\alpha +(-\Delta_v)^\alpha$ by time evolution  in
``some sense'', though the equation is only transport for $x$
variable.

We consider now a more complicate equation, the Cauchy problem for linear
Fokker-Planck equation :
\begin{equation}\label{1.00}
\left\{\begin{array}{l} f_t + v\cdot\nabla_x f =
\nabla_v\,\cdot\,\big(\nabla_v f+ v f\big), \enskip (x, v) \in
\RR^{2d},\,\,\, t
>0\, ; \\
f|_{t=0}=f_0 \, .
\end{array}\right.
\end{equation}
This equation is a natural generalization of classical Kolmogorov equation, and
a simplified model of inhomogeneous Landau equation (see \cite{villani,villani2}).
The local property of this equation is the same as classical Kolmogorov
equation since
the add terms $\nabla_v\,\cdot\,( v f )$ is a first order term,
but for the studies of kinetic equation, $v$ is velocity variable,
and hence it is in whole space $\RR^d_v$. Then
there occurs additional difficulty for analysis of this equation.

The definition of weak solution in the function space $L^\infty(]0,
T[; L^2(\RR^{2d}_{x,\,v}) \cap L^1_1(\RR^{2d}_{x,\,v}))$ for the
Cauchy problem is standard in the distribution sense, where for
$1\leq p<+\infty, l\in\RR$
$$
L^p_l(\RR^{2d}_{x,\, v})=\Big\{f\in\cS'(\RR^{2d}); (1+|v|^2)^{l/2}f\in
L^p(\RR^{2d}_{x,\, v })\Big\}.
$$
The existence of weak solution is similar to full Landau equation
(see \cite{amuxy-nonlinear2,guo}).  We get also the following
ultra-analytic effect result.

\begin{theo}\label{theo0.2} Let   $f_0\in {  L^2(\RR^{2d}_{x,\,v})
\cap L^1_1(\RR^{2d}_{x,\,v}),} 0<T\leq +\infty$. Assume that $f\in
L^\infty (]0,\,T[\-; L^2(\RR^{2d}_{x,\,v})$ $\cap
L^1_1(\RR^{2d}_{x,\,v}))$ is a weak solution of the Cauchy problem
(\ref{1.00}). Then, for any $0 < t < T$, we have
$$
f(t, \,\cdot,\,\cdot\,) \in \cA^{1/2}(\RR^{2d} ) .
$$
Furthermore, for any $0<T_0<T$ there exists $c_0>0$ such that for any
$0<t\leq T_0$, we have
\begin{equation}\label{3.4+}
\left\|e^{-c_0(t \triangle_v+t^3\triangle_x) }{f}
(t,\,\cdot,\,\cdot\,)\right\|_{L^2(\RR^{2d})} \leq  e^{ \frac {d}{2}
t} \|f_0\|_{L^2(\RR^{2d})}.
\end{equation}
\end{theo}

\begin{rema}\label{rem1-1}{
The ultra-analyticity results of the above two theorems are optimal for
the smoothness properties of solutions. From these results, we
obtain {  a good} understanding for the hypoellipticity of kinetic
equations (see \cite{desv-vil-2, hormander}), and also the
relationship, established by  Villani \cite{villani0} and
Desvillettes-Villani \cite{desv-vil-1}}, between the nonlinear Landau
equation (with Maxwellian molecules) and the linear Fokker-Planck
equation.
\end{rema}

\smallskip

We consider now the spatially inhomogeneous Landau equation
\begin{equation}\label{N-landau}
\left\{\begin{array}{l} f_t + v\cdot\nabla_x f = Q(f,\, f), \enskip
(x, v) \in \RR^{2d},\,\,\,\, t>0\, ; \\
f|_{t=0}=f_0(x, v).
\end{array}\right.
\end{equation}
The problem is now much more complicate since the solution $f$ is
the function of $(t, x, v)$ variables. We consider it here only in
the linearized framework around the normalized Maxwellian
distribution
$$
\mu (v)=(2\pi )^{-\frac{d}{2}}e^{-\frac{|v|^2}{2}},
$$
which is the equilibrium state because $Q(\mu,\, \mu)=0$.
Setting $f=\mu+g$, we consider the diffusion part of linear Landau
collision operators,
$$
Q(\mu,\,\, g)=\nabla_v \Big(\bar{a}(\mu)\cdot \nabla_v  {g}
-\bar{b}(\mu) {g}\Big)
$$
where
\begin{eqnarray*}
&&\bar{a}_{i j}(\mu)=a_{i j}\star\mu =\delta_{i j}(|v|^2+1)-v_i
v_j,\,\,\,\,\\
&&\bar{b}_j(\mu)=\sum^d_{i=1} \big(\partial_{v_i}a_{i j}\big)\star
\mu= -v_j,\,\,\,\,\,i,\,j=1,\,\cdots,\, d .
\end{eqnarray*}
In particular, it follows that
\begin{equation}
\sum^d_{i j=1}\bar{a}_{i j}(\mu) \xi_i\xi_j\geq
|\xi|^2,\,\,\,\,\mbox{for all}\,\,(v,\, \xi)\in\RR^{2d}
.\label{4.1+0}
\end{equation}
We then consider the following Cauchy problem
\begin{equation}\label{L-landau}
\left\{\begin{array}{l} g_t + v\cdot\nabla_x g = \nabla_v
\Big(\bar{a}(\mu)\cdot \nabla_v  {g} -\bar{b}(\mu) {g}\Big),\,\,
\enskip (x,\, v) \in \RR^{2d},\,\,\,\,\,\, t>0\, ; \\
g|_{t=0}=g_0.
\end{array}\right.
\end{equation}
We can also look this equation as a linear model of spatially
inhomogeneous Landau equation, which is much more complicate than
linear Fokker-Planck equation (\ref{1.00}), since the coefficients
of diffusion part are now variables. The existence and $C^\infty$
regularity of weak solution for the Cauchy problem have been considered in
\cite{amuxy-nonlinear2}. We prove now the following;

\begin{theo}\label{theo0.3}
Let $g_0\in L^2(\RR^{2d}_{x,\,v})\cap L^1_2(\RR^{2d}_{x,\,v}),
0<T\leq +\infty$. Assume that $g\in L^\infty (]0,\, T[\-;\,\,
L^2(\RR^{2d}_{x,\,v}) \cap L^1_2(\RR^{2d}_{x,\,v}) )$ is a weak
solution of  {   the } Cauchy problem (\ref{L-landau}). Then, for
any $0 < t < T$, we have
$$
g(t, \,\cdot,\,\cdot\,) \in \cA^{1}(\RR^{2d} ) .
$$
Furthermore, for any $0<T_0<T$ there {  exist} $C,\,\,c>0$ such that
for any $0<t\leq T_0$, we have
\begin{equation}\label{3.4++}
\left\|e^{c(t (-\triangle_v)^{1/2}+t^2(-\triangle_x)^{1/2}) }{g}
(t,\,\cdot,\,\cdot\,)\right\|_{L^2(\RR^{2d})} \leq e^{ C\, t}\,\,
\|g_0\|_{L^2(\RR^{2d})}.
\end{equation}
\end{theo}

In this theorem, we only consider the analytic effect result for the Cauchy
problem (\ref{L-landau}), neglecting the symmetric term
$Q(g, \mu)$ in the linearized operators of Landau collision
operator (cf.,(1.15) of \cite{amuxy-nonlinear2}) because of
the technical difficulty, see the remark in the end of section \ref{s4}.

There have been many results about the regularity of solutions
for Boltzmann equation without
angular cut-off and Landau equation, see \cite
{amuxy-nonlinear2,al-2-1, al-2-2,chen-desv-he,chen,desv-vil,desv-wen1,
HMUY,mo-xu1} for the
$C^\infty$ smoothness results, and
\cite{chen-li-xu1,chen-li-xu2,desv-fl-terr,ukai,MUXY-DCDS} for Gevrey
regularity results for Boltzmann equation and Landau equation in
both cases : the spatially homogeneous and inhomogeneous cases.
As for the analytic and Gevrey regularities, we remark that the
propagation of Gevrey regularities
of solutions is investigated in \cite{chen-li-xu2} for full nonlinear
spatially homogeneous Landau equations,
including non-Maxwellian molecule case, and the local Gevrey regularity
for all variables $t,x,v$ is
considered in \cite{chen-li-xu1}
for some semi-linear Fokker-Planck equations.
Comparing those results, the ultra-analyticity for $x,v$
variables showed in Theorem \ref{theo0.1} is
strong although the Maxwellian molecule case is only treated. As a related result
 for spatially homogeneous
Boltzmann equation in the Maxwellian molecule case, we refer
\cite{desv-fl-terr}, where the propagation of Gevrey and
ultra-analytic regularity is studied uniformly in time variable $t$.
Throughout the present paper, we focus the smoothing effect of {
the} Cauchy problem, and the uniform smoothness estimate near to
$t=0$. Concerning further details of the analytic and Gevrey
regularities of solutions  for Landau equations and Boltzmann
equation without angular cut-off, we refer the introduction of
\cite{chen-li-xu2} and references therein.

\vskip 0.5cm

\vskip0.5cm
\section{Spatially homogeneous Landau equations}\label{s2}
\setcounter{equation}{0}
\smallskip

We consider the Cauchy problem (\ref{H-landau}) and prove Theorem
\ref{theo0.1} in this section. We refer to the works of C. Villani
\cite{villani0,villani} for the essential properties of homogeneous
Landau equations. We suppose the existence of weak solution $f(t,
v)>0$ in $L^\infty(]0, T[; L^1_2(\RR^d)\-\bigcap L^2(\RR^d))$. The
conservation of mass, momentum and energy reads,
$$
\frac d {dt}\int_{\RR^d}f(t, v)\left(\begin{array}{c} 1\\ v\\ |v|^2
\end{array}\right) dv\equiv 0 .
$$
Without loss of generality, we can suppose that
\begin{eqnarray*}
&\int_{\RR^d}f(t, v) dv =1, & \mbox{unit mass}\\
&\int_{\RR^d}f(t, v) v_j dv =0,\,\,\,\,\,\,j=1, \cdots, d; &
\mbox{zero
mean velocity}\\
&\int_{\RR^d}f(t, v)|v|^2 dv =T_0,  & \mbox{unit temperature}\\
&\int_{\RR^d}f(t, v)v_j v_k dv =T_j \,\delta_{j k}, & \sum^d_j T_j=T_0\\
&T_j=\int_{\RR^d}f(t, v)v_j^2 dv >0 ,\,\,\,\,\,j=1, \cdots, d; &
\mbox{directional temperatures .}
\end{eqnarray*}
Then we have,
\begin{eqnarray}
&&\bar{a}_{j k} (f)=\delta_{j k}(|v|^2+T_0-T_j)-v_jv_k \,\,;\label{2.1}\\
&&\bar{b}_j(f)=-v_j\,\, ;\label{2.2}\\
&& \sum^d_{j, k}\bar{a}_{j k} (f)\xi_j\xi_k\geq
C_1|\xi|^2\,\,,\,\,\, \forall\,\,(v, \xi)\in \RR^{2d}\,.\label{2.3}
\end{eqnarray}
where $C_1=\min_{1\leq j\leq d}\{T_0-T_j\}>0$.

Now for $N> \frac d 4+1$ and $0<\delta<1/N,\, c_0>0,\, t>0$ , set
$$
G_\delta(t, |\xi|)=\frac{e^{c_0
t|\xi|^2}}{\left(1+\delta e^{c_0 t|\xi|^2}\right)(1+ \delta c_0 t | \xi|^2)^{N} }.
$$
Since $G_\delta(t,\, \cdot\,)\in L^\infty(\RR^d)$, we can use it as
Fourier multiplier, denoted by
$$
G_\delta(t, D_v) f (t, v)=\cF^{-1}\Big(G_\delta(t, |\xi|) \hat{f}(t,
\xi)\Big).
$$
Then, for any $t>0$,
$$
G_\delta(t)=
G_\delta(t, D_v) : L^2(\RR^d)\,\, \rightarrow\,\, H^{2N}(\RR^d)\subset C^2_b(\RR^d).
$$
The object of this section is to prove the uniform bound (with respect to
$\delta>0$) of
$$
\|G_\delta(t, D_v) f (t, \,\cdot\,)\|_{L^2(\RR^d)}.
$$

Since $f(t,\, \cdot\,)\in L^2(\RR^d) \bigcap
L^1_2(\RR^d)$ is a weak solution, we can take
$$
G_\delta(t)^2 {f} (t, \,\cdot\,)=
G_\delta(t, D_v)^2 {f} (t, \,\cdot\,)\in H^{2N}(\RR^d),
$$
as test function in the equation of (\ref{H-landau}), whence we have
\begin{eqnarray*}
&&\frac 1 2 \frac{d}{dt} \|G_\delta(t) {f}
(t,\,\cdot\,)\|^2_{L^2(\RR^d)}+ \sum^d_{j, k=1}
\int_{\RR^d}\bar{a}_{j\, k}(f)\Big(\partial_{v_j}G_\delta(t) f(t,
v)\Big)\\
&&\hskip 5cm \times\,\,
\overline{\Big(\partial_{v_k}G_\delta(t) f(t, v)\Big)}d v\\
& =&\frac 1 2\left(\Big(\partial_t G_\delta(t)\Big) f, G_\delta(t)
f\right)_{L^2(\RR^d)}+\sum^d_{j=1}
\int_{\RR^d}\Big(\partial_{v_j}\big( v_{j} f (t, v)\big)\Big)\,\,
\overline{ G_\delta(t)^2 f(t, v)}\, d v\\
&&+\sum^d_{j, k=1} \int_{\RR^d}\Big\{\bar{a}_{j
k}(f)\Big(G_\delta(t) \partial_{v_j} f(t, v)\Big)- G_\delta(t)\Big(\bar{a}_{j k}(f)
\partial_{v_j} f(t, v)\Big)\Big\}\\
&&\hskip 5cm \times\,\, \overline{\Big(\partial_{v_k}G_\delta(t)
f(t, v)\Big)}\,\, d v.
\end{eqnarray*}

To estimate the terms in the above equality, we prove the following two propositions.

\begin{prop}\label{prop2.1}
We have
\begin{eqnarray}\label{2.5}
&&C_1\|\nabla_v G_\delta(t) f(t)\|^2_{L^2(\RR^d)}\leq \sum^d_{j,
k=1} \int_{\RR^d}\bar{a}_{j k}(f)\Big(\partial_{v_j}G_\delta(t, D_v)
f(t, v)\Big)\nonumber\\
&&\hskip 5cm \times\,\, \overline{\Big(\partial_{v_k}G_\delta(t,
D_v) f(t, v)\Big)}\, d v.
\end{eqnarray}

\begin{equation}\label{2.6}
\left|\left(\Big(\partial_t G_\delta(t)\Big) f,\,\, G_\delta(t)
f\right)_{L^2}\right|\leq c_0\|\nabla_v G_\delta(t)f(t)\|^2_{L^2} .
\end{equation}

\begin{eqnarray}\label{2.6+1}
&&\re \sum^d_{j=1} \int_{\RR^d}\Big(\partial_{v_j}\big( v_{j} f (t,
v)\big)\Big)\,\, \overline{ G_\delta(t)^2 f(t, v)}\, d v \nonumber\\
&&\hskip 3cm \leq\frac d 2 \| G_\delta(t) f(t)\|^2_{L^2} + 2c_0 t
\|\nabla_v G_\delta(t) f(t)\|^2_{L^2}\, .
\end{eqnarray}
\end{prop}

\smallbreak

\noindent {\bf Proof :} The estimate
(\ref{2.5}) is exactly the elliptic condition (\ref{2.3}).
By using the Fourier transformation, (\ref{2.6}) is deduced from the
following calculus
$$
\partial_t G_\delta(t, |\xi|)=c_0|\xi|^2G_\delta(t, |\xi|)
\left(\frac {1}
 {1+\delta e^{c_0 t|\xi|^2}}-\frac {N\delta}
 {1+\delta c_0 t|\xi|^2}\right)= c_0|\xi|^2G_\delta(t, |\xi|)\, J_{N, \delta}
$$
where
$$
|J_{N, \delta}|=\left|\frac {1}
 {1+\delta e^{c_0 t|\xi|^2}}-\frac {N\delta}
 {1+\delta c_0 t|\xi|^2}\right|\leq 1.
$$
To treat (\ref{2.6+1}), we use
\begin{equation}\label{2.7}
\partial_{\xi_j} G_\delta(t, |\xi|)=2c_0 t \xi_j G_\delta(t, |\xi|)\,J_{N, \delta}.
\end{equation}
Then, we have
\begin{eqnarray*}
&&\re \sum^d_{j=1} \int_{\RR^d}\Big(\partial_{v_j}\big( v_{j} f (t,
v)\big)\Big)\,\, \overline{ G_\delta(t, D_v)^2 f(t, v)}\, d v\\
&=&-\re \sum^d_{j=1} \int_{\RR^d}v_{j} G_\delta(t, D_v) f(t, v)
\overline{\Big(\partial_{v_j} G_\delta(t, D_v) f(t, v)\Big)} \, d v\\
&&-\re \sum^d_{j=1} \int_{\RR^d}\Big([G_\delta(t, D_v), v_{j}]
f(t,v)\Big)\overline{ \Big(\partial_{v_j} G_\delta(t, D_v) f(t, v)\Big)}\, d v\\
&=&\frac d 2 \| G_\delta(t) f(t, \,\cdot)\|^2_{L^2(\RR^d)}
 -\re \sum^d_{j=1}
\int_{\RR^d}\Big([G_\delta(t, D_v), v_{j}] f(t, v)\Big)\\
&&\hskip 5cm \times\,\, \overline{\Big(\partial_{v_j} G_\delta(t,
D_v) f(t, v)\Big)}\, d v.
\end{eqnarray*}
Using Fourier transformation and (\ref{2.7}), we have that for
$t>0$,
\begin{eqnarray*}
&&-\sum^d_{j=1} \int_{\RR^3}\Big([G_\delta(t, D_v), v_{j}] f(t,
v)\Big) \overline{\Big(\partial_{v_j} G_\delta(t, D_v) f(t, v)\Big)}dv\\
&=&-\sum^d_{j=1} \int_{\RR^d}\Big(G_\delta(t, D_v) v_{j} f(t, v)
-v_j G_\delta(t, D_v) f(t, v) \Big) \overline{\Big(\partial_{v_j}
G_\delta(t, D_v) f(t, v)\Big)}dv\\
 &=&\sum^d_{j=1} \int_{\RR^d}\Big\{i
\partial_{\xi_j}\big(G_\delta(t, |\xi|) \hat{f}(t,
\xi)\Big)-G_\delta(t, |\xi|) \big(i
\partial_{\xi_j}\hat{f}(t, \xi)\big)\Big\} G_\delta(t, |\xi|)\,\,
\overline{\, i \xi_j \hat{f}(t,
\xi)}d\xi\\
&=&\sum^d_{j=1} \int_{\RR^3}\Big(\partial_{\xi_j}G_\delta(t,
|\xi|)\Big) \hat{f}(t, \xi) {\xi_j} G_\delta(t, |\xi|)
\overline{\hat{f}(t, \xi)}
d\xi\\
&=&2c_0t \int_{\RR^d}|\xi|^2 |G_\delta(t, |\xi|) \hat{f}(t, \xi)|^2
J_{N, \delta}\, d\xi
\leq 2c_0t \int_{\RR^d} |\xi|^2|G_\delta(t, |\xi|) \hat{f}(t, \xi)|^2
d\xi\,,
\end{eqnarray*}
which give (\ref{2.6+1}). The proof of Proposition
\ref{prop2.1} is now complete.

\bigbreak
For the commutator term, the special structure of the operator implies
\begin{prop}\label{prop2.2}
\begin{eqnarray*}
&&\sum^d_{j, k=1} \int_{\RR^d}\Big\{\bar{a}_{j k}(f)\Big(G_\delta(t,
D_v)
\partial_{v_j} f(t, v)\Big)- G_\delta(t, D_v)\Big(\bar{a}_{j k}(f)
\partial_{v_j} f(t, v)\Big)\Big\}\\
&&\hskip 5cm \times\,\, \overline{\Big(\partial_{v_k}G_\delta(t,
D_v) f(t, v)\Big)}\, d v=0.
\end{eqnarray*}
\end{prop}

\smallbreak

\noindent {\bf Proof :} We introduce now polar coordinates on $\RR^d_\xi$
by setting $r = |\xi|$ and
$\omega = \xi/|\xi| \in \SSS^{d-1}$. Note
that $\partial/\partial\xi_j = \omega_j \partial/
\partial r + r^{-1} \Omega_j$ where $\Omega_j$ is a vector field on
$\SSS^{d-1}$, and (see \cite{hormander}, Proposition 14.7.1)
\begin{equation}\label{2.92}
\sum_{j=1}^d \omega_j \Omega_j =0, \enskip \sum_{j=1}^d \Omega_j
\omega_j =d-1 \enskip .
\end{equation}
By using Fourier transformation, we have
\begin{eqnarray*}
&&-\sum^d_{j, k=1}\int_{\RR^d}\Big\{\bar{a}_{j k}(f)\Big(G_\delta(t,
D_v) \partial_{v_j} f(t, v)\Big)- G_\delta(t, D_v)\Big(\bar{a}_{j
k}(f)
\partial_{v_j} f(t, v)\Big)\Big\}\\
&&\hskip 6cm \times\,\, \overline{\Big(\partial_{v_k}G_\delta(t,
D_v) f(t, v)\Big)}\, d v\\
&=&\int_{\RR^d}\Big\{ \sum^d_{j, k=1} \xi_k \Big[\, \Big(\delta_{j
k}\Delta_\xi - \partial_{\xi_k}\partial_{\xi_j}\Big)\, ,\,
G_\delta(t, |\xi|) \, \Big] \, \xi_j \hat{f}(t, \xi) \Big\} \times
G_\delta(t, |\xi|) \overline{\hat{f}(t, \xi)} d\xi.
\end{eqnarray*}
Noting, in polar coordinates on $\RR^d_\xi$,
$$
\Delta_\xi = \frac{\partial^2}{\partial r^2} +
\frac{d-1}{r}\frac{\partial}{\partial r}
+ \frac{1}{r^2} \sum_{j=1}^d \Omega_j^2\,,
$$
we have, denoting by $\tilde G(r^2)=G_\delta(t,\, r)$,
\begin{eqnarray*}
&&\sum^d_{j, k=1} \omega_k \Big[\, \Big(\delta_{j k}\Big\{
\frac{\partial^2}{\partial r^2} +
\frac{d-1}{r}\frac{\partial}{\partial r}\Big\} \\
&&\hskip 2cm - \Big \{(\omega_k
\partial/ \partial r + r^{-1} \Omega_k)( \omega_j \partial/ \partial
r + r^{-1} \Omega_j)\} \Big)\, ,\, \tilde G(r^2) \, \Big] \,
\omega_j\\
&=&\Big[\, \frac{\partial^2}{\partial r^2} +
\frac{d-1}{r}\frac{\partial}{\partial r}\, ,\, \tilde G(r^2) \,
\Big] \\
&&\hskip 1cm -     \Big[\Big (\,\sum^d_{k=1} (\omega_k^2 \partial/
\partial r + r^{-1} \omega_k \Omega_k) \sum^d_{j=1}( \omega_j^2
\partial/
\partial r + r^{-1} \Omega_j \omega_j) \Big)\, ,\, \tilde G(r^2) \,
\Big]\\
&=&\Big[\, \frac{\partial^2}{\partial r^2} +
\frac{d-1}{r}\frac{\partial}{\partial r}\, ,\, \tilde G(r^2) \,
\Big] - \Big[\, \frac{\partial^2}{\partial r^2} +
\frac{\partial}{\partial r}\frac{d-1}{r}\, ,\, \tilde G(r^2) \,
\Big] =0,
\end{eqnarray*}
where we have used $(\ref{2.92})$. Then we finish the proof of
Proposition \ref{prop2.2}.

{
\begin{rema}\label{rem2}
In the above proof of Proposition \ref{prop2.2}, we have used the
polar coordinates in the dual variable of $v$, {  which} is
essentially related to a form of the Landau operator with Maxwellian
molecules. {  We notice that the same relation  (in $v$ variable)
was described by Villani \cite{villani0} and Desvillettes-Villani
\cite{desv-vil-1}}.
\end{rema}
}

\bigskip
\noindent{\bf End of proof of Theorem \ref{theo0.1}} :

\smallbreak From Propositions \ref{prop2.1} and \ref{prop2.2}, we get
\begin{eqnarray*}\displaystyle
&&\frac 1 2 \frac{d}{dt} \|G_\delta(t) {f}
(t,\,\cdot\,)\|^2_{L^2(\RR^d)}+(C_1-\frac 12 c_0-2c_0 t)\|\nabla_v
G_\delta(t) f(t,\,\cdot\,)\|^2_{L^2(\RR^d)}\\
&&\hskip 5cm \leq \frac d 2\|G_\delta(t) {f}
(t,\,\cdot\,)\|^2_{L^2(\RR^d)}.
\end{eqnarray*}
For any $0<T_0<T$,
choose $c_0$ small enough such that $C_1-\frac 12 c_0-2c_0 T_0\geq 0$.
Then 
we get
\begin{equation}\label{principle-est1+}
\frac{d}{dt} \|G_\delta(t) {f} (t,\,\cdot\,)\|_{L^2(\RR^d)}\leq
\frac d 2\|G_\delta(t) {f} (t,\,\cdot\,)\|_{L^2(\RR^d)}.
\end{equation}
Integrating the inequality (\ref{principle-est1+}) on $]0,\,t[$, we
obtain
\begin{equation}\label{2.90}
\|G_\delta(t) {f} (t,\,\cdot\,)\|_{L^2(\RR^d)}\leq e^{\frac d 2
t}\|f_0\|_{L^2(\RR^d)}.
\end{equation}
Take limit
$\delta\, \rightarrow\, 0$ in (\ref{2.90}). Then
we get
\begin{equation}\label{2.91}
\|e^{-c_0t \triangle_v} {f} (t,\,\cdot\,)\|_{L^2(\RR^d)}\leq
e^{\frac d 2 t}\|f_0\|_{L^2(\RR^d)}
\end{equation}
for any $0<t\leq T_0$. We have now proved $f(t, \,\cdot\,)\in
\cA^{1/2}(\RR^d)$ and Theorem \ref{theo0.1}.

\vskip0.5cm
\section{Linear Fokker-Planck equations}\label{s3}
\setcounter{equation}{0}
\smallskip

In the paper \cite{villani0}, there is an exact solution for
spatially homogeneous linear Fokker-Planck equation. In the
inhomogeneous case we can also obtain an exact solution of {  the}
Cauchy problem (\ref{1.00}). Denote by
$$
\hat{f}(t, \eta, \xi)=\cF_{x,
v}(f(t, x, v)),
$$
the partial Fourier transformation of $f$ with respect to $(x, v)$
variable. Then, by Fourier transformation for $(x, v)$ variables,
the linear Fokker-Planck equation (\ref{1.00}) becomes
$$\displaystyle
\left\{\begin{array}{l} \frac{\partial}{\partial t} \hat{f} (t,
\eta, \xi) -\eta\,\cdot\,\nabla_\xi \hat{f} (t, \eta, \xi)
+\xi\,\cdot\, \nabla_\xi \hat{f} (t, \eta, \xi)= -|\xi|^2\hat{f} (t,
\eta, \xi)\,; \\
 \hat{f}|_{t=0}=\cF(f_0)(\eta, \xi).
\end{array}\right.
$$
Therefore  we obtain the exact solution
$$
\hat{f}(t,\xi,\eta) = \hat{f}(0, \xi e^{-t}+ \eta(1-e^{-t}), \eta)
\,\,\exp \Big(-\int_0^t|\xi e^{\tau-t} + \eta(1-e^{\tau-t})|^2 d
\tau\Big).
$$
Note that
\begin{eqnarray*}
&&\int_0^t|\xi e^{-\tau} + \eta(1-e^{-\tau})|^2 d \tau\\
&=&\frac{1-e^{-2t}}{2}|\xi|^2 + (1-e^{-t})^2 \xi \cdot \eta + \Big(
t
- \frac{3+e^{-2t}}{2} + 2 e^{-t} \Big)|\eta|^2\\
&=&(X- \frac{X^2}{2})|\xi|^2 + X^2  \xi \cdot \eta + (-\log(1-X) -X
- \frac{X^2}{2})|\eta|^2,
\end{eqnarray*}
where $X= 1-e^{-t}\sim t$. We have for $0< K < 2/3$
$$
\int_0^t|\xi e^{-\tau} + \eta(1-e^{-\tau})|^2 d \tau \geq X(1-1/(2
K) -X/2)|\xi|^2 + (1/3 -K/2)X^3|\eta|^2.
$$
Hence for  $t \sim X < 2- 1/K$, we get
$$
f(t,\, \cdot,\,\cdot\,)\in \cA^{1/2}(\RR^{2d}),
$$
so that the ultra-analytic effect holds for any $t>0$ by means of the semi-group
property. But we cannot get the uniform estimate (\ref{3.4+}).

\smallbreak We present now  the proof of (\ref{3.4+}) which implies the
ultra-analytic effect, by commutator estimates similarly as for
homogeneous Landau equation.
Set
$$
w(t, \eta, \xi)=\hat{f}(t, \eta, \xi-t\eta).
$$
Then the Cauchy problem (\ref{1.00}) is equivalent to
\begin{equation}\label{Fokker-Planck-b}
\left\{\begin{array}{l} \frac{\partial}{\partial t} w (t, \eta, \xi)=
-|\xi-t\eta|^2 w (t, \eta, \xi)-(\xi-t\eta)\,\cdot\,\nabla_\xi
w (t, \eta, \xi)\, ; \\
w|_{t=0}=\cF(f_0)(\eta, \xi).
\end{array}\right.
\end{equation}

Since we need to study the function $\int^t_0|\xi-s\eta|^2ds$,
we prove the following
estimate.

\begin{lemm}\label{lemm4.1}
For any $\alpha>0$, there exists a constant $c_\alpha >0$ such that
\begin{equation}\label{4.1}
\int^t_0|\xi-s\eta|^\alpha ds \geq c_\alpha (t|\xi|^\alpha +
t^{\alpha+1}|\eta|^\alpha).
\end{equation}
\end{lemm}
\begin{rema} If $\alpha=2$, we can get the above
estimate by direct calculation.
The following simple proof is due to Seiji Ukai.
\end{rema}

\noindent{\bf Proof :}  Setting $s = t \tau$ and {  $\tilde \eta = t
\eta$}, we see that the estimate is equivalent to
$$
\int^1_0|\xi- \tau \tilde \eta|^\alpha d\tau \geq c_\alpha
(|\xi| ^\alpha+ |\tilde
\eta|^\alpha ).
$$
Since this is trivial when  $\tilde \eta =0$, we may assume
$\tilde \eta \ne 0$.
If $|\xi| < |\tilde \eta|$ then
\begin{eqnarray*}
&&\int^1_0|\xi- \tau \tilde \eta|^\alpha d\tau \geq |\tilde
\eta|^\alpha
\int^{1}_{0} \left| \tau - \frac{|\xi|}{|\tilde \eta|} \right|^\alpha d \tau \\
&=& |\tilde \eta|^\alpha \left \{ \int_{0}^{|\xi|/|\tilde \eta|}
\Big( \frac{|\xi|}{|\tilde \eta|}- \tau \Big)^\alpha d \tau +
\int_{|\xi|/|\tilde \eta|}^1 \Big( \tau - \frac{|\xi|}{|\tilde
\eta|} \Big)^\alpha d \tau \right\}
\\
&\geq&  \frac{|\tilde \eta|^\alpha}{\alpha+1}  \min_{0\leq \theta
\leq 1} ( \theta^{\alpha+1} + (1-\theta )^{\alpha+1})
= \frac{|\tilde \eta|^\alpha}{2^\alpha (\alpha+1)}\\
& \geq& \frac{1}{2^{\alpha+1}(\alpha+1)}(|\xi| ^\alpha+ |\tilde
\eta|^\alpha ).
\end{eqnarray*}
If $|\xi| \geq |\tilde \eta|$ then
\begin{eqnarray*}
\int^1_0|\xi- \tau \tilde \eta|^\alpha d\tau &\geq& |\xi|^\alpha
\int^{1}_0\Big( 1 - \tau \frac{|\tilde \eta|}{|\xi |}\Big)^\alpha
d\tau \geq
 |\xi|^\alpha
\int^{1}_0\big( 1 -  \tau\big)^\alpha d\tau \\
&=& \frac{|\xi|^\alpha}{\alpha+1} \geq \frac{1}{2(\alpha+1)}(|\xi|
^\alpha+ |\tilde \eta|^\alpha ).
\end{eqnarray*}
Hence we obtain (\ref{4.1}).

\smallbreak
Set now
$$
\phi(t, \eta, \xi)= c_0\left(\int^t_0|\xi-s\eta|^2ds-
\frac{c_2}{2} t^3|\eta|^2\right),
$$
where $c_0>0$ is a small constant to choose later, and $c_2$ is the constant in
(\ref{4.1}) with $\alpha=2$. Then (\ref{4.1}) implies
\begin{equation}\label{2.4}
\phi(t, \eta, \xi)\geq  c_0\frac{c_2}{2}(t|\xi|^2+ t^3|\eta|^2).
\end{equation}
Let $N =(2d+1)/4$. For $0< \delta <1/4N^2$  and $ t>0$, set
\begin{equation}\label{3.4-0}
G_\delta=G_\delta(t,\,\eta,\, \xi)= \frac{e^{\phi(t, \eta, \xi)}}{\left(1+\delta
e^{\phi(t, \eta, \xi)}\right) (1+\delta (|\eta|^2+|\xi|^2))^{ N}}\,.
\end{equation}
Since $G_\delta(t,\, \cdot\,\,\cdot\,)\in L^\infty(\RR^{2d})$, we can
use it as Fourier multiplier, denoted by
$$
\big(G_\delta(t, D_x,\, D_v) u\big) (t,\, x,\, v)=\cF^{-1}_{\eta,\,
\xi}\Big(G_\delta(t,\,\eta,\, \xi) \hat{u}(t, \,\eta,\,\xi)\Big).
$$

\begin{lemm}\label{lemm3.1}
Assume that  $f(t,\, \cdot\, ) \in L^2(\RR^{2d}_{x,\, v})
\cap L^1_1(\RR^{2d}_{x,\, v})$ for any $t\in ]0, T[$.
Then $ \nabla_\xi w(t, \,\eta,\,\xi\,)$ $\in$ $
 L^\infty(\RR^{2d}_{\eta, \xi})$, and
\begin{equation}\label{4.3+0}
|\xi-t\eta|G_\delta(t, \eta,\, \xi)^2 \bar{w} (t,
\,\eta,\,\xi\,),\,\, |\eta|G_\delta(t, \eta,\, \xi)^2 \bar{w} (t,
\,\eta,\,\xi\,),\,\, \nabla_\xi\Big( G_\delta(t, \eta,\, \xi)^2
\bar{w} (t, \,\eta,\,\xi\,)\Big)
\end{equation}
belong to $L^2(\RR^{2d}_{\eta, \xi})$ for any $t\in ]0, T[$.
\end{lemm}
\noindent{\bf Proof :}  Since $\pa_{\xi_j} w = -i \cF{( v_j f)}$, it
follows from $f \in L^1_1(\RR^{2d}_{x,\, v})$ that $ \nabla_\xi w(t,
\,\eta,\,\xi\,)$ $\in$ $
 L^\infty(\RR^{2d}_{\eta, \xi})$. Noting
$$
|\xi - t \eta|G_\delta(t, \eta,\, \xi)^2, |\eta|G_\delta(t, \eta,\, \xi)^2
\in L^\infty(\RR^{2d}_{\eta, \xi}),
$$
we see that the first two terms of (\ref{4.3+0}) are obvious.
To check the
last term in (\ref{4.3+0}), note
\begin{eqnarray}\label{2.7+}
\partial_{\xi_j} G_\delta(t, \eta, \xi)=2c_0t(\xi_j-\frac t 2
\eta_j)G_\delta(t, \eta, \xi) \frac{1}{\left(1+\delta e^{\phi(t,
\eta, \xi)}\right)}\nonumber\\
\hskip 5cm - \frac{2N \delta\xi_j}{(1+\delta
(|\eta|^2+|\xi|^2))}G_\delta(t, \eta, \xi) .
\end{eqnarray}
Then, we have
\begin{eqnarray*}
&&\nabla_\xi\Big( G_\delta(t, \eta,\, \xi)^2 \bar{w} (t,
\,\eta,\,\xi\,)\Big)=G_\delta(t, \eta,\, \xi)^2 \nabla_\xi \bar{w}
(t, \,\eta,\,\xi\,)+\nabla_\xi\Big( G_\delta(t, \eta,\, \xi)^2\Big)
\bar{w} (t,
\,\eta,\,\xi\,)\\
&&=G_\delta(t, \eta,\, \xi)^2 \nabla_\xi \bar{w} (t, \,\eta,\,\xi\,)
+4c_0t(\xi-\frac t 2 \eta) \frac{1}{\left(1+\delta
e^{\phi(t, \eta, \xi)}\right)}G_\delta(t, \eta, \xi)^2 \bar{w} (t, \,\eta,\,\xi\,)\\
&&-
\frac{ 4N \delta\xi}{(1+\delta (|\eta|^2+|\xi|^2))}
G_\delta(t, \eta, \xi)^2 \bar{w} (t, \,\eta,\,\xi\,).
\end{eqnarray*}
Since $G_\delta(t, \eta,\, \xi)^2 \in L^2(\RR^{2d}_{x,\, v})$ we have
$$
G_\delta(t, \eta,\, \xi)^2 \nabla_\xi \bar{w} (t, \,\eta,\,\xi\,)\in
L^2(\RR^{2d}).
$$
Using
$$
\left|\frac{1}{\left(1+\delta
e^{\phi(t, \eta, \xi)}\right)}\right|\leq 1,\,\,\,\,\,
\left|\frac{2N \delta\xi}{(1+\delta (|\eta|^2+|\xi|^2))}\right|\leq 1,
$$
and
\begin{eqnarray*}
&& \Big|(\xi-\frac t 2 \eta)G_\delta(t, \eta, \xi)^2 \frac{1}{\left(1+\delta
e^{\phi(t, \eta, \xi)}\right)}\bar{w} (t,
\,\eta,\,\xi\,)\Big|\leq |\xi-\frac t 2 \eta|G_\delta(t, \eta,
\xi)^2 |\bar{w} (t, \,\eta,\,\xi\,)|\\
&&\leq |\xi-t\eta|G_\delta(t, \eta, \xi)^2 |\bar{w} (t,
\,\eta,\,\xi\,)|+\frac t 2 |\eta|G_\delta(t, \eta, \xi)^2 |\bar{w}
(t, \,\eta,\,\xi\,)|\in L^2(\RR^{2d}).
\end{eqnarray*}
We have proved Lemma \ref{lemm3.1}

\bigskip

We take now $G_\delta(t, \eta,\, \xi)^2 \bar{w} (t, \,\eta,\,\xi\,)$
as test function in the equation of (\ref{Fokker-Planck-b}). Then we have
\begin{eqnarray}\label{3.0}
&&\frac{d}{dt} \|G_\delta(t, \cdot, \cdot) {w}
(t,\,\cdot\,,\, \cdot\,)\|^2_{L^2(\RR^{2d})} +
2\int_{\RR^{2d}}|(\xi-t\eta)G_\delta(t,
\eta, \xi) {w}
(t,\,\eta, \,\xi\,)|^2d\eta d\xi\nonumber\\
&&=2\sum^d_{j=1} \int_{\RR^{2d}} w (t, \eta, \xi)
\overline{\Big(\partial_{\xi_j} { (\xi_{j}-t\eta_j)} G_\delta(t,
\eta, \xi)^2 {w}(t, \eta, \xi)\Big)}d\eta d\xi\nonumber\\
&&\hskip 1cm +\left(\Big(\partial_t G_\delta(t, \cdot,\, \cdot)\Big)
w (t,\,\cdot\,,\, \cdot\,), G_\delta(t, \cdot,\, \cdot) w
(t,\,\cdot\,,\, \cdot\,)\right)_{L^2(\RR^{2d})}.
\end{eqnarray}

We prove now the following;

\begin{prop}\label{prop3.1}
We have
\begin{eqnarray}\label{3.2}
&&\left(\Big(\partial_t G_\delta(t, \cdot, \cdot)\Big) w,
G_\delta(t, \cdot, \cdot) w\right)_{L^2(\RR^{2d})} \nonumber\\
&=& c_0\int_{\RR^{2d}}|(\xi-t\eta)G_\delta(t, \eta, \xi) {w}
(t,\,\eta, \,\xi\,)|^2d\eta d\xi\nonumber\\
&&\hskip 1cm -\frac 3 2 c_0c_2 t^2\int_{\RR^{2d}}|\eta|^2
|G_\delta(t, \eta, \xi) w(t, \eta, \xi)|^2 \frac{1}{\left(1+\delta
e^{\phi(t, \eta, \xi)}\right)}\, d\eta d\xi .
\end{eqnarray}

\begin{eqnarray}\label{3.3}
&& \re \sum^d_{j=1}\int_{\RR^{2d}} w (t, \eta, \xi)
\overline{{\partial_{\xi_j}\Big((\xi_{j}-t\eta_j)}
G_\delta(t, \eta, \xi)^2  {w}(t, \eta, \xi)\Big)}d\eta d\xi \nonumber\\
&&\leq \Big( 2 c_0t+\frac{c_0 t^2}{{ 3} c_2}+c_0
\Big)\int_{\RR^{2d}}|(\xi-t\eta)G_\delta(t, \eta, \xi) {w}
(t,\,\eta, \,\xi\,)|^2d\eta d\xi \nonumber\\
&&+\frac {d +  2N^2 \delta /c_0 }{2}
\| G_\delta(t, \cdot, \cdot) w(t,\, \cdot,\,\cdot)\|^2_{L^2(\RR^{2d})}\nonumber\\
&&+\frac 3 4 c_0 c_2 t^2\int_{\RR^{2d}}|\eta|^2 |G_\delta(t, \eta,
\xi) w (t, \eta, \xi)|^2 \frac{1}{\left(1+\delta e^{\phi(t, \eta,
\xi)}\right)}\, d\eta d\xi\, .
\end{eqnarray}
\end{prop}

\noindent {\bf Proof of Proposition {\ref{prop3.1}} :} The estimate (\ref{3.2})
is deduced from 
$$
\partial_t G_\delta(t, \eta,  \xi)=c_0(|\xi-t\eta|^2-\frac 3 2 c_2t^2|\eta|^2)
G_\delta(t, \eta, \xi) \frac{1}{\left(1+\delta
e^{\phi(t, \eta, \xi)}\right)} .
$$
Since it follows from (\ref{2.7+}) that
\begin{eqnarray*}
{\mathcal I}&=&\re \sum^d_{j=1}\int_{\RR^{2d}} w (t, \eta, \xi)
\overline{{\partial_{\xi_j}\Big((\xi_{j}-t\eta_j)}
G_\delta(t, \eta, \xi)^2  {w}(t, \eta, \xi)\Big)}d\eta d\xi\\
&=&\re 2c_0t\sum^d_{j=1}\int_{\RR^{2d}}(\xi_{j}-t\eta_j)(\xi_j-\frac
t 2 \eta_j) |G_\delta(t, \eta, \xi) w (t, \eta, \xi)|^2 \frac{1}{\left(1+\delta
e^{\phi(t, \eta, \xi)}\right)}\,
d\eta d\xi\\
&&-\re \sum^d_{j=1}\int_{\RR^{2d}}
\frac{2N \delta\xi_j(\xi_{j}-t\eta_j)}{(1+\delta (|\eta|^2+|\xi|^2))}
 |G_\delta(t, \eta, \xi) w (t, \eta, \xi)|^2\,
d\eta d\xi\\
&&-\re \sum^d_{j=1}\int_{\RR^{2d}}(\xi_{j}-t\eta_j)\Big(\partial_{\xi_j}
G_\delta(t, \eta, \xi)w (t, \eta, \xi)\Big)
\overline{
G_\delta(t, \eta, \xi) {w}(t, \eta, \xi)}d\eta d\xi,
\end{eqnarray*}
we get
\begin{eqnarray*}
{\mathcal I}
&=& 2 c_0t\sum^d_{j=1}\int_{\RR^{2d}}(\xi_{j}-t\eta_j)(\xi_j-\frac t
2 \eta_j) |G_\delta(t, \eta, \xi) w (t, \eta, \xi)|^2 \frac{1}{\left(1+\delta
e^{\phi(t, \eta, \xi)}\right)}\,
d\eta d\xi\\
&-& \sum^d_{j=1}\int_{\RR^{2d}}
\frac{2 N \delta\xi_j(\xi_{j}-t\eta_j)}{(1+\delta (|\eta|^2+|\xi|^2))}
 |G_\delta(t, \eta, \xi) w (t, \eta, \xi)|^2\,
d\eta d\xi\\
&&\hskip 6cm +\frac d 2 \| G_\delta(t, \cdot, \cdot)
w(t,\,\cdot,\,\cdot)\|^2_{L^2(\RR^{2d})}
\\
&=& 2 c_0t\int_{\RR^{2d}}|\xi-t\eta|^2 |G_\delta(t, \eta, \xi) w (t,
\eta, \xi)|^2 \frac{1}{\left(1+\delta
e^{\phi(t, \eta, \xi)}\right)}\,
d\eta d\xi\\
&+&c_0t^2\int_{\RR^{2d}}(\xi-t\eta)\,\cdot\,\eta |G_\delta(t, \eta,
\xi) w (t, \eta, \xi)|^2\frac{1}{\left(1+\delta
e^{\phi(t, \eta, \xi)}\right)}\, d\eta d\xi\\
&-& \sum^d_{j=1}\int_{\RR^{2d}}
\frac{2 N \delta\xi_j(\xi_{j}-t\eta_j)}{(1+\delta (|\eta|^2+|\xi|^2))}
 |G_\delta(t, \eta, \xi) w (t, \eta, \xi)|^2\,
d\eta d\xi\\
&&\hskip 6cm +\frac d 2 \| G_\delta(t, \cdot, \cdot)
w(t,\,\cdot,\,\cdot)\|^2_{L^2(\RR^{2d})}.
\end{eqnarray*}
 For the last term, noting
$$\sum^d_{j=1} \frac{2 N \delta \xi_j(\xi_{j}-t\eta_j)}
{(1+\delta (|\eta|^2+|\xi|^2))}
\leq \frac{(N^2 /c_0)\delta^2 |\xi|^2 + c_0|\xi -t \eta|^2}
{(1+\delta (|\eta|^2+|\xi|^2))}
\leq N^2 \delta/c_0 + c_0|\xi -t \eta|^2$$
we finally obtain
\begin{eqnarray*}
{\mathcal I}
&\leq& \Big(2 c_0t+\frac{c_0 t^2}{{3} c_2} + c_0 \Big)
\int_{\RR^{2d}}|(\xi-t\eta)G_\delta(t,
\eta, \xi) {w}
(t,\,\eta, \,\xi\,)|^2d\eta d\xi\\
&&+\frac {d+2N^2 \delta/ c_0 } {2}
\| G_\delta(t, \cdot, \cdot) w(t, \cdot, \cdot)\|^2_{L^2(\RR^{2d})}\\
&&+\frac 3 4 c_0 c_2 t^2\int_{\RR^{2d}}|\eta|^2 |G_\delta(t, \eta, \xi) w
(t, \eta, \xi)|^2 \frac{1}{\left(1+\delta
e^{\phi(t, \eta, \xi)}\right)}\, d\eta
d\xi\, .
\end{eqnarray*}
Thus we have proved Proposition \ref{prop3.1}.

\bigskip
\noindent{\bf End of proof of Theorem \ref{theo0.2}} :

Now the equation (\ref{3.0}), the estimate (\ref{3.2}) and
(\ref{3.3}) deduce
\begin{eqnarray*}
&&\frac{d}{dt} \|G_\delta(t, \cdot, \cdot) {w}
(t,\,\cdot\,\cdot\,)\|^2_{L^2(\RR^{2d})}\\
&&\hskip 1cm + \Big(2-3 c_0-4 c_0t-\frac{{2}c_0 t^2}{{3} c_2}\Big)
\int_{\RR^{2d}}|(\xi-t\eta)G_\delta(t, \eta, \xi) {w}
(t,\,\eta, \,\xi\,)|^2d\eta d\xi\\
&\leq&  (d+2N^2 \delta/c_0)\,  \| G_\delta(t, \cdot, \cdot) w(t,
\cdot, \cdot)\|^2_{L^2(\RR^{2d})}.
\end{eqnarray*}
Then for any $0<T_0<T$ choose $c_0>0$
(depends on $T_0$) small enough such that
$$
2-3 c_0-4c_0 T_0-\frac{{ 2}c_0 T^2_0}{{ 3} c_2}\geq 0,
$$
then for any $0<t\leq T_0$,
$$
\frac{d}{dt} \|G_\delta(t, \cdot, \cdot) {w}
(t,\,\cdot,\,\cdot\,)\|_{L^2(\RR^{2d})}\leq  \frac {d+2N^2 \delta/c_0}{2}
\| G_\delta(t, \cdot, \cdot) w(t,\,\cdot,\,\cdot\,)\|_{L^2(\RR^{2d})},
$$
which gives
$$
\|G_\delta(t, \cdot, \cdot) {w}
(t,\,\cdot,\,\cdot\,)\|_{L^2(\RR^{2d})}\leq e^{ \frac {d +2N^2\delta/c_0}{2}t}
\|f_0\|_{L^2(\RR^{2d})}.
$$
Take $\delta\rightarrow0$, we have
\begin{eqnarray*}
&&\int_{\RR^{2d}}e^{c_0\int^t_0|\xi-s\eta|^2ds -c_1t^3|\eta|^2} |\hat{f}
(t,\,\eta,\,\xi-t\eta\,)|^2 d\eta d\xi\\
&=&\int_{\RR^{2d}}e^{c_0\int^t_0|\xi+(t-s)\eta|^2ds -c_1t^3|\eta|^2}
|\hat{f} (t,\,\eta,\,\xi\,)|^2 d\eta d\xi \leq  e^{ {d \, t}}
\|f_0\|^2_{L^2(\RR^{2d})}.
\end{eqnarray*}
By using (\ref{2.4}), we get finally
$$
\|e^{-\tilde{c}_0(t \triangle_v+t^3\triangle_x )}{f}
(t,\,\cdot,\,\cdot\,)\|_{L^2(\RR^{2d})} \leq  e^{ \frac {d}{2} t}
\|f_0\|_{L^2(\RR^{2d})}
$$
for any $0<t\leq T_0$, where $\tilde{c}_0=\frac{c_0c_2}{2}>0$.
This is the desired estimate
(\ref{3.4+}), which implies
$$
f(t,\,\cdot,\,\cdot\,)\in {\mathcal A}^{1/2}(\RR^{2d}).
$$
We have
thus proved Theorem \ref{theo0.2}.

\vskip0.5cm
\section{Linear model of inhomogeneous Landau equations}\label{s4}
\setcounter{equation}{0}
\smallskip

We prove now the Theorem \ref{theo0.3} in this section.  By the
change of variables $(t,x,v) \longrightarrow (t,\, x+v\, t,\, v)$, the
Cauchy problem (\ref{L-landau}) is reduced to
\begin{equation}\label{IH-landau}
\left\{\begin{array}{l} {f}_t=(\nabla_v - t
\nabla_x)\Big(\bar{a}(\mu)\cdot(\nabla_v
-t \nabla_x) {f} -\bar{b}(\mu) {f}\Big),\\
{f}|_{t=0}=g_0(x, v),\end{array}\right.
\end{equation}
where ${f}(t,\, x,\, v)=g(t,\, x+v\, t,\, v)$.  Recall that
\begin{eqnarray*}
&&\bar{a}_{i j}(\mu)=a_{i j}\star\mu =\delta_{i j}(|v|^2+1)-v_i
v_j\,\,;\\
&& \bar{b}_j(\mu)=\sum^d_{i=1} \big(\partial_{v_i}a_{i j}\big)\star
\mu= -v_j\,\,;\,\,\,\,\,\,i,\,j=1, \cdots, d ,
\end{eqnarray*}
and
$$
\sum^d_{i j=1}\bar{a}_{i j}(\mu) \xi_i\xi_j\geq
|\xi|^2,\,\,\,\,\,\,\,\mbox{for all}\,\,\,(v,\, \xi)\in\RR^{2d} .
$$
In view of this Cauchy problem , we set
$$
\Psi(t, \eta, \xi)= c_0 \int^t_0|\xi-s\eta|ds ,
$$
for a sufficiently small $c_0>0$ which will be chosen later on. Then
we can use the (\ref{4.1}) with $\alpha=1$ to estimate $\Psi$. Set
$$
F_\delta(t,\,\eta,\, \xi)= \frac{e^{\Psi}}{(1+\delta e^{\Psi})(1+
\delta \Psi)^N }
$$
for $N=d+1, 0<\delta \leq \frac 1{N}$. 
If $A$ is a first
order differential operator of $(t,\eta,\xi)$ variables then we have
\begin{equation}\label{4.3}
A F_{\delta} =
 \left( \frac{1}{1+\delta e^{\Psi}}
 -\frac{N\delta}{1+ \delta \Psi}\right)
 (A \Psi) F_{\delta},
\end{equation}
and
$$
\left| \frac{1}{1+\delta e^{\Psi}}
 -\frac{N\delta }{1+ \delta \Psi}\right|\leq 1 .
$$

\smallskip
Taking
$$
F_\delta(t,D_x, D_v)^2 f= F_\delta(t)^2 f \in H^{2N}(\RR^{2d})
$$
as a test
function in the weak solution formula of (\ref{IH-landau}), we have
\begin{eqnarray*}
&&\frac 1 2 \frac{d}{dt} \|F_\delta(t) {f}
\|^2_{L^2(\RR^{2d})}+\left(\bar{a}(\mu)\Big((\nabla_{v}- t
\nabla_{x})F_\delta(t) f\Big),\, \Big((\nabla_{v}- t\nabla_{x})
F_\delta(t) f\Big)\right)_{L^2(\RR^{2d})}\\
&=&\displaystyle -\sum^d_{j=1} \int_{\RR^{2d}}v_{j} f\,\,
\overline{\Big((
\partial_{v_j}  - t \partial_{x_j})F_\delta(t)^2 f  \Big)}dx dv
+\frac 1 2\Big(\big(\partial_t F_\delta \big) f, \, F_\delta(t)
f\Big)_{L^2(\RR^{2d})}\\
&& +\displaystyle \sum^d_{j, k=1} \int_{\RR^{2d}}\Big\{\bar{a}_{j
k}(\mu)\Big(F_\delta(t) (\partial_{v_j} - t \partial_{x_j} )
\Big)\,f- F_\delta(t)\Big(\bar{a}_{j k}(\mu) (\partial_{v_j} - t
\partial_{x_j}) f\Big)\Big\} \\
&&\hskip 5cm \times \,\,\overline{\Big((\partial_{v_k} - t
\partial_{x_k}) F_\delta(t) f \Big)}dx dv.
\end{eqnarray*}

We prove now the following results.

\begin{prop}\label{prop4.1}
We have
\begin{eqnarray}
&& \|(\nabla_v - t \nabla_x) F_\delta(t) f
\|^2_{L^2(\RR^{2d})}\nonumber \\
&&\hskip 1cm
\leq\left(\bar{a}(\mu)\Big((\nabla_{v}- t \nabla_{x})F_\delta(t)
f\Big),\, \Big((\nabla_{v}- t\nabla_{x}) F_\delta(t)
f\Big)\right)_{L^2(\RR^{2d})} .\label{4.5}\\
&&\label{4.6} \left|\left(\Big(\partial_t F_\delta(t) \Big) f,\,
F_\delta(t)
 f\right)_{L^2}\right|\leq c_0 \|(\nabla_v - t \nabla_x) F_\delta(t)
f\|_{L^2} \|F_\delta(t) f\|_{L^2}. \\
&&\label{4.6+1} -\re \sum^d_{j=1} \int_{\RR^{2d}}v_{j} f\,\,\,
\overline{\Big((
\partial_{v_j}  - t \partial_{x_j})F_\delta(t)^2 f  \Big)}
\leq\frac d 2 \| F_\delta(t) f\|^2_{L^2}  \hskip3cm  \nonumber\\
&&\hskip5cm + c_0 t \|(\nabla_v - t \nabla_x) F_\delta f(t)\|_{L^2}
\|F_\delta f(t) \|_{L^2} \, .
\end{eqnarray}
\end{prop}

\smallbreak
\noindent {\bf Proof :} The estimate (\ref{4.5}) is a direct consequence of the
elliptic condition (\ref{4.1+0}). Using the Fourier transformation
and noting (\ref{4.3}), we see that (\ref{4.6})  is derived {}from
$$
\partial_t F_{\delta}(t,\eta,\xi) =
 \left( \frac{1}{1+\delta e^{\Psi}}
 -\frac{N\delta }{1+ \delta \Psi}\right)
 (\partial_t \Psi) F_{\delta}\,, \enskip \partial_t \Psi = c_0 |\xi - t
 \eta|.
$$

For (\ref{4.6+1}), we have firstly
$$
-\re \sum^d_{j=1} \int_{\RR^6}v_{j} F_\delta(t) f\,\,\,
\overline{\Big((
\partial_{v_j}  - t \partial_{x_j})F_\delta(t) f  \Big)}
= \frac d 2 \| F_\delta(t) f\|^2_{L^2}.
$$
For the commutators $[v_{j},\,  F_\delta(t)]$, using Fourier
transformation, we have that for $t>0$ and $\hat f = \hat
f(t,\eta,\xi)$
\begin{eqnarray*}
&&-\sum^d_{j=1} \int_{\RR^{2d}}\Big([F_\delta(t,D_x,D_v), v_{j}]
f(t,x, v)\Big) \overline{\Big((\partial_{v_j}- t\partial_{x_j})
F_\delta(t, D_x, D_v) f(t,x, v)\Big)}dx dv\\
&&=-\sum^d_{j=1} \int_{\RR^{2d}}\Big(F_\delta(t, D_x, D_v) v_{j}
f(t) -v_j F_\delta(t, D_x, D_v) f(t) \Big)\\
&&\hskip 5cm \,\times\,\, \overline{\Big((\partial_{v_j}- t
\partial_{x_j})
F_\delta(t, D_v) f(t)\Big)}dx dv\\
 &&=\sum^3_{j=1} \int_{\RR^{2d}}\Big\{i
\partial_{\xi_j}\big(F_\delta(t, \eta, \xi)
\hat{f}(t)\Big)-F_\delta(t, \eta, \xi) \big(i
\partial_{\xi_j}\hat{f}(t)\big)\Big\} F_\delta(t, \eta, \xi)\,\,
\\
&&\hskip 7cm \,\times\,\, \overline{\, i (\xi_j - t \eta_j) \hat{f}(t)} d\eta d\xi\\
&&=\sum^d_{j=1} \int_{\RR^{2d}}\Big(\partial_{\xi_j}F_\delta(t,
\eta, \xi)\Big) \hat{f}(t) {(\xi_j- t \eta_j)} F_\delta(t, \eta,
\xi) \overline{\hat{f}(t)}
d \eta d\xi\\
&&\leq  c_0 t\int_{\RR^{2d}}|\xi -t \eta| |F_\delta(t, \eta, \xi)
\hat{f}(t)|^2 d\eta d\xi\leq c_0 t\|(\nabla_v - t \nabla_x) F_\delta
f(t)\|_{L^2} \|F_\delta f(t) \|_{L^2} ,
\end{eqnarray*}
where, in view of (\ref{4.3}), we have used the fact that
$$
\left |\sum_{j=1}^d(\partial_{\xi_j} \Psi) (t,\eta, \xi) \times
(\xi_j - t \eta_j)\right| \leq c_0 \int_0^1 \left| \sum_{j=1}^3
\frac{\xi_j - s \eta_j}{|\xi - s \eta| } (\xi_j - t \eta_j)\right|
ds \leq c_0 t |\xi - t \eta|.
$$
Thus (\ref{4.6+1}) has been proved.

\bigbreak
For the commutator terms, we have

\begin{prop}\label{prop4.2}
There exists a constant $C_1 >0$ independent of $\delta>0$ such that
\begin{eqnarray}\label{4.6+2}
&&\Big| \sum^d_{j, k=1} \int_{\RR^{2d}}\Big\{\bar{a}_{j k}(\mu)
\Big(F_\delta(t) (\partial_{v_j} - t \partial_{x_j} ) \Big)\,f-
F_\delta(t)\Big(\bar{a}_{j k}(\mu) (\partial_{v_j} - t
\partial_{x_j}) f\Big)\Big\}\nonumber\\
 && \hskip 6cm \times\,\,\, \overline{\Big((\partial_{v_k} - t \partial_{x_k})
F_\delta(t) f \Big)}\,\,\, \Big|\nonumber\\
&&\leq C_1 \left\{
 (c_0 t)^2\|(\nabla_v - t
\nabla_x) F_\delta(t) f\|_{L^2} ^2 + \|F_\delta(t) f\|_{L^2} ^2
\right\}.
\end{eqnarray}
\end{prop}

\smallbreak
\noindent {\bf Proof  :} In order to prove (\ref{4.6+2}), we introduce the polar
coordinates of $\xi$ centered at $t\eta$, that is ,
$$
r = |\xi -t \eta|\,\,\,\,\,\, \mbox{and}\,\,\,\,\,\, \omega =
\frac{\xi-t\eta}{|\xi - t \eta|} \in \SSS^{d-1}.
$$
Note again  that $\partial/\partial\xi_j = \omega_j \partial/
\partial r + r^{-1} \Omega_j$ where $\Omega_j$ is a vector field on
$\SSS^{d-1}$. We have again
$$
\sum_{j=1}^d \omega_j \Omega_j =0, \enskip \sum_{j=1}^d \Omega_j
\omega_j =d-1, \enskip .
$$
By means of Plancherel formula, we have
\begin{eqnarray*}
&&\displaystyle \sum^d_{j, k=1} \int_{\RR^{2d}}\Big\{\bar{a}_{j
k}(\mu) \Big(F_\delta(t) (\partial_{v_j} - t \partial_{x_j} ) \Big)-
F_\delta(t)\Big(\bar{a}_{j k}(\mu) (\partial_{v_j} - t
\partial_{x_j}) f\Big)\Big\}\\
&&\hskip 7cm \times\,\, \overline{\Big((\partial_{v_k}
- t \partial_{x_k}) F_\delta(t) f \Big)}\\
&=&- \int_{\RR^{2d}}\Big\{ \sum^d_{j, k=1} (\xi_k- t \eta_k) \Big[\,
\Big(\delta_{j k}\Delta_\xi -
\partial_{\xi_k}\partial_{\xi_j}\Big)\, ,\, F_\delta (t, \eta, \xi)\, \Big] \,
(\xi_j - t \eta_j) \hat f(t) \Big\} \\
&&\hskip 7cm \times\,\,   \overline{ F_\delta(t,
\eta, \xi) \hat f(t)} d\xi d \eta \\
&=& J\,.
\end{eqnarray*}
Noting again
$$
\Delta_\xi = \frac{\partial^2}{\partial r^2} +
\frac{d-1}{r}\frac{\partial}{\partial r}
+ \frac{1}{r^2} \sum_{l=1}^d \Omega_l^2\,,
$$
we have with $\tilde F_\delta(t, \eta,\, r,\, \omega )= F_\delta(t,
\eta, r\cdot \omega + t\eta)=F_\delta (t, \eta, \xi)$
\begin{eqnarray*}
&&-\sum^d_{j, k=1} \omega_k \left[\, \Big(\delta_{j k}\Big\{
\frac{\partial^2}{\partial r^2} +
\frac{d-1}{r}\frac{\partial}{\partial r} + \frac{1}{r^2}
\sum_{l=1}^d \Omega_l^2\Big\}\right.\\
&&\hskip 2cm \,\,  -\left.\Big \{(\omega_k \frac{\partial}{
\partial r }+ r^{-1} \Omega_k)( \omega_j \frac{\partial}{ \partial r
}+ r^{-1} \Omega_j) \Big \} \Big)\, ,\,\,\, \tilde F_{\delta} \,
\right]
\, \omega_j\\
&&=-\Big[\, \frac{\partial^2}{\partial r^2} +
\frac{d-1}{r}\frac{\partial}{\partial r}\, ,\,\,\, \tilde F_\delta
\, \Big]\\
&&\hskip 2cm \,\,+\Big[\Big (\,\sum^d_{k=1} (\omega_k^2
\frac{\partial}{ \partial r} + r^{-1} \omega_k \Omega_k)
\sum^d_{j=1}( \omega_j^2 \frac{\partial}{
\partial r }+ r^{-1} \Omega_j \omega_j) \Big)\, ,\, \,\,\tilde F_\delta \,
\Big] \\
&&-\frac{1}{r^2} \sum^d_{j=1} \omega_j \Big[\,  \sum_{l=1}^d
\Omega_l^2 ,\,\,\, \tilde F_{\delta} \, \Big] \, \omega_j = A_1 +
A_2 +A_3.
\end{eqnarray*}
Note again that
$$
A_1 + A_2 =-\Big[\, \frac{\partial^2}{\partial r^2} +
\frac{d-1}{r}\frac{\partial}{\partial r}\, ,\, \tilde F_\delta\,
\Big] +\Big[\, \frac{\partial^2}{\partial r^2} +
\frac{\partial}{\partial r}\frac{d-1}{r}\, ,\, \tilde F_\delta \,
\Big] =0.
$$
On the other hand, we have in view of (\ref{4.3})
\begin{eqnarray*}
&&A_3 = -\frac{1}{r^2}\sum^d_{j, l=1} \omega_j \Big(2 \Omega_l
[\Omega_l ,\,\, \tilde F_\delta] - \Big[\,  \Omega_l\,, \,\,
[\Omega_l ,\, \,\tilde F_{\delta}]\Big]  \Big) \,
\omega_j \\
&&= -\frac{1}{r^2}\sum^d_{j, l=1}
 \omega_j \left ( 2 \Omega_l
\left(\frac{(\Omega_l \Psi) }{1+ \delta e^{\Psi}} - \frac{N \delta
(\Omega_l \Psi) }{
1+ \delta \Psi} \right) \tilde F_\delta \right. \\
&&\left. - \left( \left( \frac{(\Omega_l \Psi) }{1+ \delta e^{\Psi}}
- \frac{N \delta  (\Omega_l \Psi) } {1+ \delta \Psi} \right)^2 +
\left(
 \Omega_l \left(\frac{(\Omega_l \Psi) }{1+ \delta e^{\Psi}} - \frac{N \delta
(\Omega_l \Psi) } {1+ \delta \Psi} \right) \right) \right) \tilde
F_\delta \right ) \omega_j.
\end{eqnarray*}
Putting $w_j =  \omega_j \tilde F_\delta\,\, w$ with $w(t, \eta, r,
\omega) =\hat{f}(t, \eta, r\,\cdot\,\omega+t\eta)$, we have
\begin{eqnarray*}
&J &= \re J = \re \int_{\RR^d_\eta} \int_0^\infty \int_{S^{d-1}}
r^2 (A_3 w)\,\, \overline{ \tilde F_\delta w}\,\, r^{d-1} d r d \omega d \eta\\
&=&- \sum^d_{j, l=1} \mbox{Re} \int_{\RR^d_\eta} \int_0^\infty
\int_{S^{d-1}} \left\{2 \Omega_l \left(\frac{(\Omega_l \Psi) }{1+
\delta e^{\Psi}} - \frac{N \delta (\Omega_l \Psi) }{ 1+ \delta \Psi}
\right)
w_j \right\}  \overline{w_j}\,\, r^{d-1} d r d \omega d \eta\\
&+& \sum^d_{j, l=1}\int_{\RR^d_\eta} \int_0^\infty \int_{S^{d-1}}
\left( \left( \frac{(\Omega_l \Psi) }{1+ \delta e^{\Psi}} - \frac{N
\delta  (\Omega_l \Psi) } {1+ \delta \Psi} \right)^2\right.\\
&& \hskip 3cm +\left. \left(
 \Omega_l
 \left(\frac{(\Omega_l \Psi) }{1+ \delta e^{\Psi}} - \frac{N \delta
(\Omega_l \Psi) } {1+ \delta \Psi} \right) \right) \right)
|w_j|^2 r^{d-1} d r d \omega d \eta\\
&=& J_1 + J_2.
\end{eqnarray*}
Since $\Omega^*_l  = - \Omega_l + (d-1)\omega_l$,
the integration by parts gives
\begin{eqnarray*}
J_1 = & -\sum^d_{j, l=1}\int_{\RR^d_\eta} \int_0^\infty \int_{S^{d-1}}
\left\{ \left(
 \Omega_l \left(\frac{(\Omega_l \Psi) }{1+ \delta e^{\Psi}} - \frac{N \delta
(\Omega_l \Psi) } {1+ \delta \Psi} \right) \right)  \right. \\
&+ (d-1) \omega_l
\left.
\left(\frac{(\Omega_l \Psi) }{1+ \delta e^{\Psi}} - \frac{N \delta
(\Omega_l \Psi) } {1+ \delta \Psi} \right) \right\}
|w_j|^2 r^{d-1} d
r d \omega d \eta.
\end{eqnarray*}
Hence we obtain
\begin{eqnarray}\label{Y-4.1}
J&=& \sum^d_{j, l=1}\int_{\RR^d_\eta} \int_0^\infty \int_{S^{d-1}}
\left\{ \left( \frac{1}{1+ \delta e^{\Psi}} - \frac{N \delta  } {1+ \delta
\Psi} \right)^2 (\Omega_l \Psi) ^2 \right.\\
&&- (d-1) \omega_l
\left.
\left(\frac{1}{1+ \delta e^{\Psi}} - \frac{N \delta
 } {1+ \delta \Psi} \right)(\Omega_l \Psi) \right\}
|w_j|^2 r^{d-1} d
r d \omega d \eta \nonumber \\
&=&\int_{\RR^d_\eta} \int_0^\infty \int_{S^{d-1}} \left\{
\left( \frac{1}{1+
\delta e^{\Psi}} - \frac{N \delta  } {1+ \delta \Psi} \right)^2
\left(\sum^d_{l=1}(\Omega_l \Psi)^2\right) \right.\nonumber \\
&&\left. -(d-1)
\left( \frac{1}{1+
\delta e^{\Psi}} - \frac{N \delta  } {1+ \delta \Psi} \right)
\left(\sum^d_{l=1}\omega_l (\Omega_l \Psi)\right)\right\}
|\tilde F_\delta w |^2
r^{d-1} d r d \omega d \eta. \nonumber
\end{eqnarray}
Since there exists a constant $C_d >0$ such that
\begin{equation}\label{Y-4.2}
|\Omega_l \Psi| =
c_0 r \left|\sum_{j=1}^d \int_0^t \frac{\xi_j - s \eta_j}{|\xi - s \eta|} ds
(\Omega_l \omega_j) \right| \leq c_0 C_d t r ,
\end{equation}
we have
\begin{eqnarray*}
|J|&\leq&  C_d'\left \{(c_0 t)^2 \int_{\RR^d_\eta} \int_0^\infty
\int_{S^{d-1}} r^2 |\tilde F_\delta w |^2 r^{d-1} d r d \omega d
\eta\right.\\
&&\hskip 3cm + \left.\int_{\RR^d_\eta} \int_0^\infty \int_{S^{d-1}}
|\tilde F_\delta w |^2 r^{d-1} d r d \omega d \eta\right\} \,,
\end{eqnarray*}
which yields  (\ref{4.6+2}). The proof of Proposition \ref{prop4.2}
is now complete.

\bigskip
\noindent {\bf End of proof of Theorem \ref{theo0.3} : }

\smallskip

{}From Propositions \ref{prop4.1} and \ref{prop4.2}, there
exist constants $C_2,\, C_3>0$ independent of $\delta>0$ and
$t> 0$ such that
\begin{eqnarray*}
&&\frac 1 2 \frac{d}{dt} \|(F_\delta{f})
(t)\|^2_{L^2(\RR^{2d})}+\left(\frac 1 2-(c_0 t)^2 C_2
\right)\|(\nabla_v - t \nabla_x) (F_\delta
f)(t)\|^2_{L^2(\RR^{2d})}\\
&&\hskip 6cm \leq C_3 \|(F_\delta {f} )(t)\|^2_{L^2(\RR^{2d})}.
\end{eqnarray*}
So that if $\frac 1 2 -(c_0 t)^2 C_2 \geq 0$, we have,
\begin{equation}\label{principle-est1}
\frac{d}{dt} \|(F_\delta {f}) (t)\|_{L^2(\RR^{2d})}\leq
 C_3 \|(F_\delta{f} )(t)\|_{L^2(\RR^{2d})}.
\end{equation}
Using the fact $(F_\delta {f}) (0)=\frac{1}{1+\delta}\, g_0$, we get
$$
\|(F_\delta{f} )(t)\|_{L^2(\RR^{2d})}\leq e^{C_3
t}\|g_0\|_{L^2(\RR^{2d})}.
$$
Take the limit $\delta\, \rightarrow\, 0$. Then we have
\begin{equation}\label{4.a11}
\int_{\RR^{{2d}}}e^{2\Psi(t, \eta, \xi)}|\hat{f} (t, \eta,\xi)|^2
d\eta d\xi\leq e^{2C_3 t}\|g_0\|^2_{L^2(\RR^{2d})}.
\end{equation}
On the other hand, by Lemma \ref{lemm4.1}, there exists a $c_1
>0$ such that
\begin{eqnarray*}
\int_{\RR^{2d}}e^{2\Psi(t, \eta, \xi)}|\hat{f} (t, \eta,\xi)|^2
d\eta d\xi&=&\int_{\RR^{2d}}e^{2c_0\int^t_0|\xi-s\eta|ds}|\hat{g}
(t,
\eta, \xi-t\eta)|^2 d\eta d\xi\\
&=&\int_{\RR^{2d}}e^{2c_0\int^t_0|\xi+(t-s)\eta|ds}|\hat{g} (t,
\eta,
\xi)|^2 d\eta d\xi\\
&\geq&\int_{\RR^{2d}}e^{2c_0 c_1 (t|\xi| + t^2|\eta|)}|\hat{g} (t,
\eta, \xi)|^2 d\eta d\xi .
\end{eqnarray*}
Finally, for any $0<T_0<T$, choosing $c_0>0$ small enough such that $\frac
1 2 -(c_0 T_0)^2 C_2 \geq 0$, we have proved,
$$
\int_{\RR^{2d}}\left|e^{c_0 c_1 (t(-\triangle_v)^{1/2} +
t^2(-\triangle_x)^{1/2})}g (t,\, x,\, v)\right|^2 dx\, dv \leq
e^{2C_3 t}\|g_0\|^2_{L^2(\RR^{2d})}\,\,\,\,\,\, \mbox{for any}\,\,\,
0<t\leq T_0\,,
$$
which completes the proof of Theorem \ref{theo0.3} with $C=2 C_3$
depending only on $d$.

\begin{rema}\label{rem4}
The formulas (\ref{Y-4.1}) and (\ref{Y-4.2}) show that we cannot
get the ultra-analytic effect of order $1/2$ as in Theorem \ref{theo0.2}. It is
the same reason why we do not consider the symmetric term $Q(g, \mu)$
in the equation (\ref{L-landau}) as in \cite{amuxy-nonlinear2}.
\end{rema}

\bigskip\noindent
{\bf Acknowledgments:} Authors wish to express their hearty gratitude
to Seiji Ukai who communicated Lemma \ref{lemm4.1}.
The research of the first author was
supported by  Grant-in-Aid for Scientific Research No.18540213,
Japan Society of the Promotion of Science, and the
second author would like to thank the  support of
Kyoto University for his visit there.

\end{document}